\documentclass[10pt,a4paper]{article}
\usepackage{cite}
\usepackage{amscd}
\usepackage[latin1]{inputenc}
\usepackage{amsmath}
\usepackage{amsfonts}
\usepackage{amssymb}
\usepackage{color}
\usepackage{float}
\usepackage{amsthm}
\usepackage{graphicx}
\usepackage{listings}
\newtheorem{theorem}{Theorem}

\newtheorem{proposition}[theorem]{Proposition}

\def\m{\mathbb}
\def\mc{\mathcal}

\begin{document}
\date{}
\title{A surface with $q=2$ and canonical\\ map of degree $16$}
\author{Carlos Rito}
\maketitle

\begin{abstract}

We construct a surface with irregularity $q=2,$ geometric genus $p_g=3,$
self-intersection of the canonical divisor $K^2=16$ and canonical map of
degree $16.$

\noindent 2010 MSC: 14J29.

\end{abstract}

\section{Introduction}

Let $S$ be a smooth minimal surface of general type.
Denote by ${\phi:S\dashrightarrow\m P^{p_g-1}}$ the canonical map and let $d:=\deg(\phi).$
The following Beauville's result is well-known.
\begin{theorem}[\cite{Be}]
If the canonical image $\Sigma:=\phi(S)$ is a surface, then either:
\begin{description}
\item{{\rm (i)}} $p_g(\Sigma)=0,$ or
\item{{\rm (ii)}} $\Sigma$ is a canonical surface $($in particular $p_g(\Sigma)=p_g(S))$.
\end{description}
Moreover, in case {\rm (i)} $d\leq 36$ and in case {\rm (ii)} $d\leq 9.$
\end{theorem}

Beauville has also constructed families of examples with $\chi(\mathcal O_S)$ arbitrarily
large for $d=2, 4, 6, 8$ and $p_g(\Sigma)=0.$ Despite being a classical problem, for $d>8$ the
number of known examples drops drastically: only Tan's example \cite[\S 5]{Ta} with
$d=9,$ the author's \cite{Ri} example with $d=12$ and Persson's example \cite{Pe} with $d=16$ are known.
Du and Gao \cite{DuGa} show that if the canonical map is an abelian cover
of $\m P^2,$ then these examples with $d=9$ and $d=16$ are the only possibilities for $d>8.$
These three surfaces are regular, so for irregular surfaces all known examples satisfy $d\leq 8$.
We get from Beauville's proof that lower bounds hold for irregular surfaces.
In particular, $$q=2\ \Longrightarrow\ d\leq 18.$$

In this note we construct an example with $q=2$ and $d=16.$ The idea of the construction is the following.
We start with a double plane with geometric genus $p_g=3,$ irregularity $q=0,$
self-intersection of the canonical divisor $K^2=2$ and singular set the union of $10$ points of type
$\sf A_1$ (nodes) and $8$ points of type $\sf A_3$ (standard notation, the resolution
of a singularity of type $\sf A_n$ is a chain of $(-2)$-curves $C_1,\ldots,C_n$ such that $C_iC_{i+1}=1$ and
$C_iC_j=0$ for $j\not=i\pm 1$). Then we take a double covering ramified over the points of type $\sf A_3$
and obtain a surface with $p_g=3,$ $q=0$ and $K^2=4$ with $28$ nodes.
A double covering ramified over $16$ of these $28$ nodes gives a surface with $p_g=3,$ $q=0$ and $K^2=8$ with
$24$ nodes (which is a $\m Z_2^3$-covering of $\m P^2$).
Finally there is a double covering ramified over these $24$ nodes which gives
a surface with $p_g=3,$ $q=2$ and $K^2=16$ and the canonical map factors through these coverings,
thus it is of degree $16.$

\bigskip
\noindent{\bf Notation}

We work over the complex numbers. All varieties are assumed to be projective algebraic.
A $(-n)$-curve on a surface is a curve isomorphic to $\m P^1$ with self-intersection $-n.$
Linear equivalence of divisors is denoted by $\equiv.$
The rest of the notation is standard in Algebraic Geometry.\\

\bigskip
\noindent{\bf Acknowledgements}

The author thanks Thomas Baier for many interesting conversations.

The author is a member of the Center for Mathematics of the University of Porto.
This research was partially supported by FCT (Portugal) under the project PTDC/MAT-GEO/0675/2012 and by CMUP (UID/MAT/00144/2013), which is funded by FCT with national (MEC) and European structural funds through the programs FEDER, under the partnership agreement PT2020.

\section{$\m Z_2^n$-coverings}\label{coverings}

The following is taken from \cite{Ca}, an alternative reference is \cite{Pa}.

\begin{proposition}
A normal finite $G\cong\m Z_2^r$-covering $Y\rightarrow X$ of a smooth variety $X$ is completely
determined by the datum of
\begin{enumerate}
\item reduced effective divisors $D_{\sigma},$ $\forall\sigma\in G,$ with no common components;
\item divisor classes $L_1,\ldots,L_r,$ for $\chi_1,\ldots,\chi_r$ a basis of the dual group of characters
$G^{\vee},$ such that $$2L_i\equiv\sum_{\chi_i(\sigma)=-1}D_{\sigma}.$$ 
\end{enumerate}
Conversely, given 1. and 2., one obtains a normal scheme $Y$ with a finite
$G\cong\m Z_2^r$-covering $Y\rightarrow X.$
\end{proposition}
The covering $Y\rightarrow X$ is embedded in the total space of the direct sum of the line bundles whose
sheaves of sections are the $\mc O_X(-L_i),$ and is there defined by equations
$$u_{\chi_i}u_{\chi_j}=u_{\chi_i+\chi_j}\prod_{\chi_i(\sigma)=\chi_j(\sigma)=-1}x_{\sigma},$$
where $x_{\sigma}$ is a section such that ${\rm div}(x_{\sigma})=D_{\sigma}.$
The scheme $Y$ can be seen as the normalization of the Galois covering given by the equations
$$u_{\chi_i}^2=\prod_{\chi_i(\sigma)=-1}x_{\sigma}.$$
The scheme $Y$ is irreducible if $\{\sigma|D_{\sigma}>0\}$ generates $G.$

For the reader's convenience, we leave here the character table for the group $\m Z_2^3$ with generators $x,y,z.$
{\footnotesize
\begin{verbatim}
[    -1    -1    -1    -1     1     1     1     1 x*y*z]
[    -1    -1     1     1    -1    -1     1     1     z]
[    -1     1    -1     1    -1     1    -1     1     y]
[    -1     1     1    -1     1    -1    -1     1     x]
[     1    -1    -1     1     1    -1    -1     1   y*z]
[     1    -1     1    -1    -1     1    -1     1   x*z]
[     1     1    -1    -1    -1    -1     1     1   x*y]
[     1     1     1     1     1     1     1     1    Id]
\end{verbatim}
}

\section{The construction}

\noindent{\bf Step 1}\\
Let $T_1,\ldots,T_4\subset\m P^2$ be distinct lines tangent to a smooth conic $H_1$
and $$\pi:X\longrightarrow\m P^2$$ be the double cover of the projective plane ramified over $T_1+\cdots+T_4.$
The curve $\pi^*(H_1)$ is of arithmetic genus $3,$ from the Hurwitz formula, and has $4$ nodes,
corresponding to the tangencies to $T_1+\cdots+T_4$. Hence $\pi^*(H_1)$ is reducible,
$$\pi^*(H_1)=A+B$$ with $A,B$ smooth rational curves. From $AB=4$ and $(A+B)^2=8$ we get $A^2=B^2=0.$
Now the adjunction formula $$2g(A)-2=AK_X+A^2$$ gives $AK_X=-2$ and then the Riemann-Roch Theorem implies
$$h^0(X,\mc O_X(A))\geq 1+\frac{1}{2}A(A-K_X)=2.$$ 
Therefore there exists a smooth rational curve $C$ such that $C\not=A,$ $C\equiv A$ and $AC=0.$
The curve $$H_2:=\pi(C)$$ is smooth rational. The fact $\pi^*(H_2)^2>C^2$ implies that $\pi^*(H_2)$ is reducible, 
thus $H_2$ is tangent to the lines $T_1,\ldots,T_4.$ As above, there is a smooth rational curve $D$
such that $$\pi^*(H_2)=C+D$$ and $C^2=D^2=0.$
Since $A\equiv C$ and $A+B\equiv C+D,$ then $B\equiv D.$\\
\begin{figure}[h]
\centering
\includegraphics[width=6.5cm,height=3.3cm]{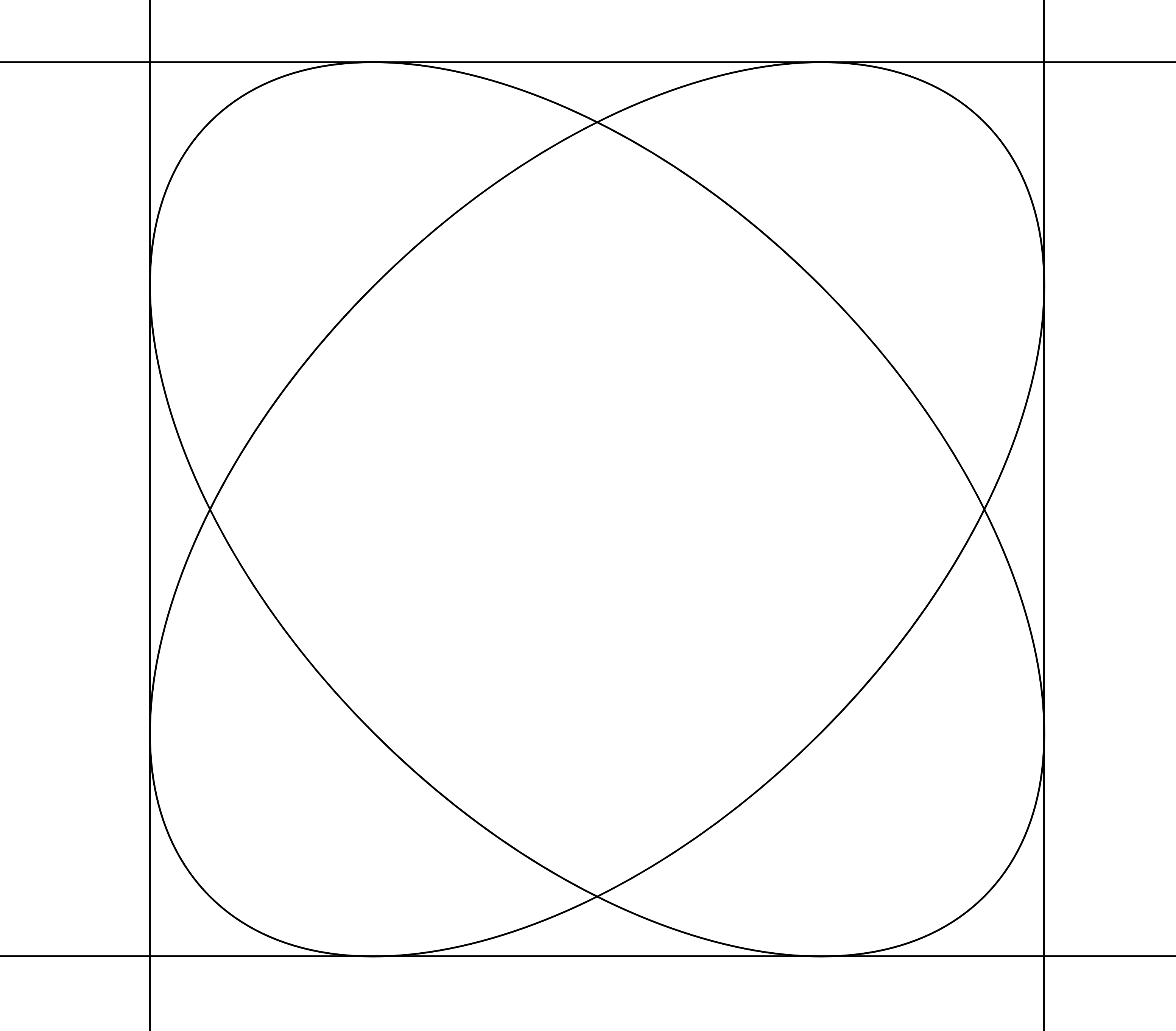}
\vspace*{-0.5cm}
$$H_1+H_2+T_1+\cdots+T_4$$
\end{figure}

\noindent{\bf Step 2}\\
Let $x, y, z$ be generators of the group $\m Z_2^3$ and $$\psi:Y\longrightarrow\m P^2$$ be the
$\m Z_2^3$-covering defined by
$$D_1:=D_{xyz}:=H_1,\ D_2:=D_z:=H_2,\ D_3:=D_y:=T_1+T_2,\ D_4:=D_x:=T_3+T_4,$$
$$D_{yz}:=D_{xz}:=D_{xy}:=0.$$
Let $d_i$ be the defining equation of $D_i.$
According to Section \ref{coverings}, the surface $Y$ is obtained as the normalization
of the covering given by equations
$$u_1^2=d_1d_2d_3d_4,\ u_2^2=d_1d_2,\ \ldots,\ u_7^2=d_3d_4.$$

Since the branch curve $D_1+\cdots+D_4$ has only negligible singularities, the invariants of $Y$ can
be computed directly. Consider divisors $L_{i\ldots h}$ such that $2L_{i\ldots h}\equiv D_i+\cdots+D_h$
and let $T$ be a general line in $\mathbb P^2.$ We have
$$L_{1234}(K_{\m P^2}+L_{1234})=4T\cdot T=4,$$
$$L_{ij}(K_{\m P^2}+L_{ij})=2T(-T)=-2,$$
thus
$$\chi(Y)=8\chi\left(\m P^2\right)+\frac{1}{2}\left(4+6\times(-2)\right)=4,$$
$$p_g(Y)=p_g\left(\m P^2\right)+h^0\left(\m P^2,\mc O_{\m P^2}(T)\right)+6h^0\left(\m P^2,\mc O_{\m P^2}(-T)\right)=3.$$
So a canonical curve in $Y$ is the pullback of a line in $\m P^2$ and then $$K_Y^2=8.$$\\

\noindent{\bf Step 3}\\
Notice that the points where two curves $D_i$ meet transversely give rise to smooth points of $Y$, hence the
singularities of $Y$ are:
\begin{description}
\item[$\cdot$] $16$ points $p_1,\ldots,p_{16}$ corresponding to the tacnodes of $D_1+\cdots+D_4;$
\item[$\cdot$] $8$ nodes $p_{17},\ldots,p_{24}$ corresponding to the nodes of $D_3$ and $D_4.$
\end{description}
We want to show that $p_1,\ldots,p_{24}$ are nodes with even sum.

The surface $X$ defined in Step 1 is the double plane with equation $u_7^2=d_3d_4,$
thus the covering $\psi$ factors trough a $\m Z_2^2$-covering $$\varphi:Y\longrightarrow X.$$
The branch locus of $\varphi$ is $A+B+C+D$ plus the $4$ nodes given by the points in $D_3\cap D_4.$
The points $p_1,\ldots,p_{16}$ are nodes because they are the pullback of nodes of $A+B+C+D$.

The divisor $\varphi^*(A+C)$ is even ($A+C\equiv 2A$), double ($A+C$ in the branch locus of $\varphi$),
with smooth support ($A+C$ smooth) and $p_1,\ldots,p_{16}\in \varphi^*(A+C),$ $p_{17},\ldots,p_{24}\notin \varphi^*(A+C).$
Consider the minimal resolution of the singularities of $Y$
$$\rho:Y'\longrightarrow Y$$ and let $A_1,\ldots,A_{24}\subset Y'$ be the $(-2)$-curves
corresponding to the nodes $p_1,\ldots,$ $p_{24}.$
The divisor $(\varphi\circ\rho)^*(A+C)$ is even and there exists a divisor $E$ such that
$$(\varphi\circ\rho)^*(A+C)=2E+\sum_1^{16}A_i.$$ Thus there exists a divisor $L_1$ such that
$\sum_1^{16}A_i\equiv 2L_1.$

Analogously one shows that the nodes $p_{17,}\ldots,p_{24}$ have even sum, i.e. there exists a divisor $L_2$
such that $\sum_{17}^{24}A_i\equiv 2L_2.$ This follows from $\psi^*(T_1+T_3)$ even, double, and with support
of multiplicity $1$ at $p_{17,}\ldots,p_{24}$ and of multiplicity $2$ at $8$ of the nodes $p_1,\ldots,p_{16}.$\\

\noindent{\bf Step 4}\\
So there is a divisor $L:=L_1+L_2$ such that $$\sum_{1}^{24}A_i\equiv 2L.$$
Consider the double covering $S\longrightarrow Y$ ramified over $p_1,\ldots,p_{24}$ and determined by $L.$
More precisely, given the double covering $$\eta:S'\longrightarrow Y'$$ with branch locus $\sum_{1}^{24}A_i,$
determined by $L,$ $S$ is the minimal model of $S'$.
We have $$\chi(S')=2\chi(Y')+\frac{1}{2}L(K_{Y'}+L)=8-6=2.$$

Since the canonical system of $Y$ is given by the pullback of the system of lines in $\m P^2,$  
the canonical map of $Y$ is of degree $8$ onto $\m P^2.$ We want to show that the canonical map
of $S'$ factors through $\eta.$

One has $$p_g(S')=p_g(Y')+h^0(Y',\mc O_{Y'}(K_{Y'}+L)),$$ so the canonical map factors if 
$$h^0(Y',\mc O_{Y'}(K_{Y'}+L))=0.$$ Let us suppose the opposite. Hence the linear system $|K_{Y'}+L|$
is not empty and then $A_i(K_{Y'}+L)=-1,$ $i=1,\ldots,24,$ implies that $\sum_1^{24}A_i\equiv 2L$
is a fixed component of $|K_{Y'}+L|.$ Therefore
$$h^0(Y',\mc O_{Y'}(K_{Y'}+L-2L))=h^0(Y',\mc O_{Y'}(K_{Y'}-L))>0$$
and then
$$h^0\left(Y',\mc O_{Y'}\left(2K_{Y'}-2L\right)\right)=
h^0\left(Y',\mc O_{Y'}\left(2K_{Y'}-\sum_1^{24}A_i\right)\right)>0.$$
This means that there is a bicanonical curve $B$ through the $24$ nodes of $Y.$
We claim that there is exactly one such curve.
In fact, the strict transform in $Y'$ of the line $T_1$ is the union of two double curves $2T_a, 2T_b$
such that $$T_a\sum_1^{24}A_i=T_b\sum_1^{24}A_i=6$$ and $T_a\rho^*(B)=T_b\rho^*(B)=4.$
This implies that $\rho^*(B)$ contains $T_a$ and $T_b.$
Analogously $\rho^*(B)$ contains the reduced strict transform of $T_2, T_3$ and $T_4.$
There is only one bicanonical curve with this property, with equation $u_7=0$
(the bicanonical system of $Y$ is induced by $\mc O_{\m P^2}(2)$ and $u_2,\ldots,u_7$).

As $$h^0(Y',\mc O_{Y'}(2K_{Y'}-2L))=1\ \Longrightarrow\ h^0(Y',\mc O_{Y'}(K_{Y'}-L))=1,$$
then such bicanonical curve is double. This is a contradiction because the curve given by $u_7=0$
is not double.

So $h^0(Y',\mc O_{Y'}(K_{Y'}+L))=0$ and we conclude that the surface $S$ has invariants $p_g=3,$ $q=2,$
$K^2=16$ and the canonical map of $S$ is of degree $16$ onto $\m P^2.$

\bibliography{ReferencesRito}

\

\

\noindent Carlos Rito\\
\\{\it Permanent address:}
\\ Universidade de Tr\'as-os-Montes e Alto Douro, UTAD
\\ Quinta de Prados
\\ 5000-801 Vila Real, Portugal
\\ www.utad.pt
\\ crito@utad.pt\\
\\{\it Current address:}
\\ Departamento de Matem\' atica
\\ Faculdade de Ci\^encias da Universidade do Porto
\\ Rua do Campo Alegre 687
\\ 4169-007 Porto, Portugal
\\ www.fc.up.pt
\\ crito@fc.up.pt

\end{document}